\documentclass[a4paper,12pt,reqno]{article}

\usepackage{amsthm}
\usepackage{amsmath,amssymb,latexsym,amsfonts,mathrsfs}
\usepackage[dvipdfmx]{graphicx}
\usepackage{bm}

\usepackage{color}
\usepackage{comment}
\usepackage{mathtools}
\usepackage{fancyhdr}
\usepackage[top=30truemm,bottom=30truemm,left=20truemm,right=20truemm]{geometry}
\usepackage{enumerate}

\newcommand{\R}{\mathbb{R}}

\newcommand{\e}{\varepsilon}

\newcommand{\SCR}[1]{{\mathscr #1}}

\newcommand{\CAL}[1]{{\cal #1}}

\theoremstyle{plain}
\newtheorem{Thm}{{\bf Theorem}}[section]

\newtheorem{Lem}[Thm]{{\bf Lemma}}
\newtheorem{Prop}[Thm]{{\bf Proposition}}

\newtheorem{Ass}[Thm]{{\bf Assumption}}

\theoremstyle{definition}
\newtheorem{Def}[Thm]{{\bf Definition}}
\newtheorem{Rem}[Thm]{{\bf Remark}}

\newcounter{Exami}

\allowdisplaybreaks[2]
\begin{document}
\fontencoding{T1}\selectfont
\title{Short time asymptotics of the fundamental solutions for Schr\"{o}dinger equations with non-smooth potentials}
\author{Shun Takizawa}
\date{\today}
\maketitle
\begin{abstract}
This paper deals with Schr\"{o}dinger equations with potentials which are time-dependent non-smooth and at most quadratic growth. In the case where potentials are smooth with respect to spatial variables, fundamental solutions have explicit formulas in short time by D. Fujiwara. On the otherhand in the case where ones are non-smooth, we cannot expect that fundamental solutions have similar formula as above because dispersive estimates fail to hold in general. We show that a principal part of an asymptotic form of the fundamental solution has similar form as above even in the case where a potential is in $C^2$ with respect to spatial variables.
\end{abstract}
\section{Introduction} 
We consider short time asymptotic behaviors of the fundamental solution of the following Cauchy problem for Schr\"{o}dinger equations with non-smooth and unbounded potentials:
%
\begin{align}\label{CP1}
\begin{cases}
&i\partial_{t}u=-\frac{1}{2}\Delta u+V(t,x)u, \hspace{3mm}(t,x)\in \R \times \R^n,
\\&u(0,x)=u_{0}(x),\hspace{3mm}x\in\R^n,
\end{cases}
\end{align}
where $i=\sqrt{-1}$, $u(t,x)$ is a complex-valued function of $(t,x)\in \R \times \R^n$, $u_0(x)$ is a complex-valued function of $x\in\R^n$, $V(t,x)$ is a real-valued function, $\partial_{t}=\partial/\partial t$ and $\Delta=\sum_{j=1}^{n}\partial^2/\partial x_{j}^2$.
%
%
%
\begin{Ass}\label{AssV}
The potential $V(t,x)$ satisfies two conditions below:
\begin{itemize}
\item (Regularity) For every $t\in\R$ fixed, $V(t,x)$ belongs to $C^{2}(\R^n)$ with respect to $x$. 
Moreover $\partial_{x}^{\alpha}V(t,x)$ is continuous with respect to $t$ and $x$ for multi-indices $\alpha$ with $|\alpha|\leq2$.
\item (Growth rate) For multi-indices $\alpha$ with $|\alpha|=2$, there exists $C_{\alpha}>0$ such that
\begin{equation*}
|\partial_{x}^{\alpha}V(t,x)|\leq C_{\alpha},\hspace{3mm}(t,x)\in \R\times\R^n.
\end{equation*}
\end{itemize}
\end{Ass}
Under Assumption \ref{AssV}, there exists a unique solution of initial value problem \eqref{CP1} in $C(\R, L^2(\R^n))$. Hence the propagator $U(t)$ uniquely exists (see \cite{Nicola} Proposition 4.2).
We call the distribution kernel of $U(t)$ fundamental solution, which we write $E(t,\cdot,\cdot)$ for:
\begin{equation*}
u(t,x)=U(t)u_0(x)=\int_{\R^n}E(t,x,y)u_0(y)dy.
\end{equation*}

In \cite{Fujiwara1, Fujiwara2}, D. Fujiwara has shown that if $V$ is smooth in addition to Assumption \ref{AssV}, then the fundamental solution $E(t,x,y)$ has the form in short time:
\begin{equation}\label{3}
E(t,x,y)=(2\pi it)^{-n/2}e^{i S(t,x,y)}a(t,x,y),
\end{equation}
where $a(t,x,y)=1+O(t),\hspace{2mm}0<t\ll1$ in $L_{x,y}^{\infty}$ and $S$ is the action integral defined by 
\begin{equation*}
S(t,x,y)=\int_{0}^{t}\left\{\frac{1}{2}|\xi(s)|^2-V(s,x(s))\right\}ds,
\end{equation*}
where $x(s)$ and $\xi(s)$ are solutions to the Hamiltonian system corresponding to (\ref{CP1}):
\begin{equation*}
\dot{x}(s)=\xi(s),\hspace{3mm}\dot{\xi}(s)=-\left(\nabla_{x}V\right)\left(s,x(s)\right)
\end{equation*}
with boundary conditions $x(0)=y, \hspace{2mm}x(t)=x$.

The formula (\ref{3}) gives rigorous mathematical meaning Feynman path integral in \cite{Feynman1, Feynman2}.
In short time, the formula (\ref{3}) is a generalization of the well-known formulae for the fundamental solutions for without potential ($V=0$), with Stark potential ($V=Ex$) and with the harmonic potential ($V=x^2/2$).
However we cannot expect that the formula (\ref{3}) holds for non-smooth potentials $V$. 
It is because dispersive estimate: $\|U(t)\|_{L^1\to L^{\infty}}\lesssim t^{-n/2}, 0<t\ll1$ may fail to hold in the case time independent potentials $V(x)\in C^{\alpha}(\R^n), \alpha<(n-3)/2$. where $V$ are compactly supported and $n>3$ (see M. Goldberg-M. Visan \cite{Goldberg}).
The purpose of this paper is to show that a principal part of short time asymptotic form for the fundamental solution has similar form to (\ref{3}) and to examine a convergence rate of remainder term of one when $V$ is in $C^2$ with respect to spatial variables. 

%
%
%
%
The following theorem is our main result.
\begin{Thm}\label{main thm}
Under Assumption \ref{AssV}, we let $E(t,x,y)$ be the fundamental solution for (\ref{CP1}). Then there exist $T>0$ and positive constant $C_1, C_2$ such that 
\begin{equation}
E(t,x,y)=(2\pi it)^{-\frac{n}{2}}e^{i S(t,x,y)}a_0(t,x,y)+r(t,x,y), \hspace{3mm}(t,x,y)\in(0,T)\times\R^n\times\R^n,
\label{main-form}
\end{equation}
\begin{equation}\label{inqa_0}
\left\|a_0(t,x,y)-1\right\|_{L_{x,y}^{\infty}}\leq C_1 t,\hspace{3mm}0<t<T,
\end{equation}
\begin{equation}\label{inqr}
\left\| \int_{\R^n} r(t,\cdot,y)f(y)dy\right\|_{L^{2}(\R^n)} \leq C_2  t^2 \|f\|_{L^{2}(\R^n)},\hspace{3mm}0<t<T,\hspace{3mm}f\in L^{2}(\R^n).
\end{equation}
\end{Thm}
%

%
%
\begin{Rem}
We remark on conditions in Assumption \ref{AssV} and Theorem \ref{main thm}.
The condition of short time $t\in(0,T)$ is essential in order to guarantee well-definedness of the action integral $S$.
We cannot relax the assumption on growth rate of $V$, that is, Theorem \ref{main thm} fails to hold in the case where $V$ is superquadratic. It is because behaviors of the fundamental solutions dramatically change when the growth rate at infinity of $V$ changes from quadratic to superquadratic 
(see K. Yajima \cite{Yajima1, Yajima2}).
\end{Rem}

\begin{Rem}
F. Nicola \cite{Nicola} has proved that the fundamental solution has the form  \eqref{main-form} with $a_0\equiv1$ in the case that $\partial_{x}^{\alpha}V(t,\cdot)\in L^{\infty}(\R; H_{ul}^{n+1}(\R^n))$ with $|\alpha|=2$, where $H_{ul}^{n+1}(\R^n)$ is uniformly local $L^2$-Sobolev space. In \cite{Nicola}, the $L^2$ boundedness theorem for oscillatory integral operators by A. Boulkhemair \cite{Boulk} is used. Hence the condition regarding regularity of  $V$ cannot be relaxed in the similar method to \cite{Nicola}.
\end{Rem}
\vskip\baselineskip
\textbf{Notation.}
Let $\nabla_x=(\partial_{x_1},\ldots,\partial_{x_n})$.
We use notation $\SCR{F}f(\xi)=\widehat{f}(\xi)=\int_{\R^n}f(x)e^{-ix\cdot\xi}dx$ for Fourier transform of $f$ 
and $\SCR{F}^{*}f(x)=\int_{\R^n}f(\xi)e^{ix\cdot\xi}\bar{d}\xi$ with $\bar{d}\xi=(2\pi)^{-n}d\xi$ for the  inverse Fourier transform of $f$.
We define the dilation $D_{\e}$ which leaves $L^2$-norm invariant by $\left(D_{\e}f\right)(x)=\e^{-n/2}f(\e^{-1}x), \hspace{2mm}\e>0$.
We write $X\lesssim Y$ if $X\leq CY$ with some harmless constant $C>0$ in the proofs.
We often write $L^p$ instead of $L^p(\R^n)$ or $L^p(\R^{2n})$ for short and write $xy$ for the Euclid inner product $x\cdot y$.
We define the Schr\"{o}dinger propagator of a free particle $e^{\frac{1}{2}it\Delta}$ by 
\begin{align} 
e^{\frac{1}{2}it\Delta}f(x)=\SCR{F}^{*}e^{-\frac{1}{2}it\xi^2}\SCR{F}f(x) 
=\int_{\R^n}(2\pi it)^{-n/2}e^{i\frac{(x-y)^2}{2t}}f(y)dy. \label{9}
\end{align}

\section{Preliminaries}
%
%
%
\subsection{Wave packet transform} 
We recall the definition of the wave packet transform due to C\'{o}rdoba-Fefferman \cite{CF}.
\begin{Def}
Let $\varphi \in \CAL{S}(\mathbb{R}^n)\setminus \{0\}$ and $f \in \CAL{S} '(\mathbb{R}^n)$. Then the wave packet transform $W_{\varphi}f$ of $f$ with respect to $\varphi$ is defined by
\begin{equation*}
W_{\varphi}f(x, \xi)=\int_{{\mathbb R}^n} \overline{\varphi(y-x)}f(y)e^{-iy\cdot \xi}\bar{d}y,\quad (x, \xi) \in{\mathbb R}^n\times{\mathbb R}^n.
\end{equation*}
We also define the formal adjoint operator $W_{\varphi}^{*}$ of $W_{\varphi}$ by
\begin{equation*}
W_{\varphi}^{*}F(x)=\iint_{{\mathbb R}^{2n}} \varphi(x-y)F(y, \xi)e^{ix\cdot \xi}dy\bar{d}\xi,\quad x\in{\mathbb R}^n, F\in \CAL{S}^{'}(\R^{2n}),
\end{equation*}
with $\bar{d}\xi=(2\pi)^{-n}d\xi$. 
\end{Def}
We briefly recall properties of wave packet transform in order to prove Theorem \ref{main thm}.
\begin{Prop}\label{Prop}
Let $\varphi \in \CAL{S}(\mathbb{R}^n)\setminus \{0\}$. Then it follows that
\begin{enumerate}
\item $f(x)=\|\varphi\|_{L^2}^{-2} W_{\varphi}^{*}[W_{\varphi}f](x)$, $f\in{\CAL S}'({\R}^n)$.
%
%
\item $\|W_{\varphi}^{*}F\|_{L^{2}(\R^n)}\leq (2\pi)^{-n}\|\varphi\|_{L^{2}(\R^n)} \|F\|_{L^{2}(\R^{2n})}$, $F\in L^{2}(\R^{2n})$.
\end{enumerate}
\end{Prop}
\begin{proof}
It is known that first statement is valid by (\cite{Grochenig} Corollary 11.2.7).
Concerning the second statement, we consider the case $F\in \CAL{S}(\mathbb{R}^{2n})$ by density. The Schwarz inequality yields
\begin{equation*}
|W_{\varphi}^{*}F(x)|\leq\|\varphi\|_{L^{2}(\R^n)} \|\SCR{F}^{*}F(\cdot,x)\|_{L^{2}(\R^n)}.
\end{equation*}
Taking $L_{x}^{2}(\R^n)$-norm of the both sides of this inequality and using Plancherel's theorem, we have the desired inequality.
\end{proof}
%
%
%
%
\subsection{Classical orbits}
We consider the following Hamiltonian system corresponding to (\ref{CP1}):
\begin{equation}\label{H}
\dot{x}(s)=\xi(s), \hspace{2mm}\dot{\xi}(s)=-(\nabla_{x}V)(s,x(s)).
\end{equation}
We write $x^{(1)}(s;t,x,y),\xi^{(1)}(s;t,x,y)$ for solutions to (\ref{H}) with the boundary condition $x(0)=y, x(t)=x$ and
write $x^{(2)}(s;t,x,\xi),\xi^{(2)}(s;t,x,\xi)$ for solutions to (\ref{H}) with the initial condition $x(t)=x, \xi(t)=\xi$.
We often omit the paramaters and write $x^{(1)}(s)$ etc.
If the potential $V$ satisfies Assumption \ref{AssV}, then there exists a unique global solution $\left(x^{(2)}(s), \xi^{(2)}(s)\right)$.
Otherwhile, under Assumption \ref{AssV}, there exists a unique solution $\left(x^{(1)}(s), \xi^{(1)}(s)\right)$
in short time: $t<\pi/\sqrt{M}$, where $M=\sup\{|\partial_{x}^{2}V(t,x)|; (t,x)\in \R^{1+n}, |\alpha|=2\}$
by Narita-Ozawa \cite{Narita-Ozawa}. 
%
%
%
%
\begin{Lem}
\label{Lem1}  
Under Assumption \ref{AssV}, there exist $T>0$ and $r_0\in  \SCR{B}((0,T)\times \R^n \times \R^n)$ such that 
\begin{equation*}
\left|\frac{\partial (\xi)}{\partial (x^{(2)}(0;t,x,\xi))}\right|=t^{-n}+t^{-n+1}r_0(t,x,X)
\end{equation*}
for $0<t<T$, where 
$\frac{\partial (\xi)}{\partial (x^{(2)}(0;t,x,\xi))}$ is the Jacobian of $x^{(2)}(0;t,x,\xi)$ as the function $\xi \longmapsto X=x^{(2)}(0;t,x,\xi)$.
\end{Lem}
\begin{proof}
There exist $T_1>0$ and $C>0$ such that $|\partial_{\xi_j}x^{(2)}_k(0)+\delta_{j,k}t|\leq C t^2 $ for $0<t<T_1$ (see Yajima \cite{Yajimabook} Lemma 7.14.).
Thus we have $\partial_{\xi_j}x^{(2)}_k(0)=-t\delta_{j,k}+t^2 r_{j,k}(t,x,\xi)$, where $r_{j,k}$ is bounded.
Hence we obtain 
\begin{equation*}
\frac{\partial (x^{(2)}(0;t,x,\xi))}{\partial (\xi)}=(-t)^{n}+t^{n+1}r'_0(t,x,\xi),
\end{equation*}
where $r'_0$ is bounded. We put $M'=\sup\{|r'_0(t,x,\xi)|; (t,x,\xi)\in (0,T)\times \R^{n+n}\}$ and $T=\min\{T_1, (2M')^{-1}\}$.
Let $t\in (0,T)$ fixed. There exists constant $C>0$ independent of $x, \xi$ such that
\begin{equation*}
\left|\frac{\partial (x^{(2)}(0;t,x,\xi))}{\partial (\xi)}\right|>C.
\end{equation*}
By Hadamard's inverse mapping theorem, 
\begin{equation*}
\frac{\partial (\xi)}{\partial (x^{(2)}(0;t,x,\xi))}=\{(-t)^{n}+t^{n+1}r'_0(t,x,X)\}^{-1}=(-t)^{-n}+t^{-n+1}\widetilde{r_0}(t,x,X),
\end{equation*}
with $\widetilde{r_0}=\frac{(-1)^{n+1}r'_0}{1+(-1)^ntr'_0}$.
Let 
$r_0(t,x,X)=t^{n-1}\left|\frac{\partial(\xi)}{\partial (x^{(2)}(0))}\right|-t^{-1}$.
Then $|r_0|\leq |\widetilde{r_0}|$ for all $(t,x,X)\in (0,T)\times \R^{n+n}$.
\end{proof}
%
%
\begin{Lem} \label{Lem2}
Under Assumption \ref{AssV}, let $0<T<\pi/\sqrt{M}$ fixed. Then there exist $C_1,C_2>0$ such that
\begin{align} 
&|x^{(1)}(s;t,x,y)-x^{(2)}(s;t,x',\xi)|\leq C_1 \left(|x-x'|+|y-x^{(2)}(0;t,x',\xi)|\right), \label{inq1}
\\&
\left|\xi^{(1)}(s;t,x,y)-\xi^{(2)}(s;t,x',\xi)-\frac{(x-x')-(y-x^{(2)}(0;t,x',\xi))}{t}\right| \notag
\\&\leq C_2 t \left(|x-x'|+|y-x^{(2)}(0;t,x',\xi)|\right) \label{inq2},
\end{align}
 for all $0\leq s\leq t, 0<t<T$ and $x, y, x', \xi\in\R^n$.
\end{Lem}
\begin{proof}
Let $T<\min\{1,\pi/\sqrt{M}\}$ fixed. Integrating the both sides of the (\ref{H}), we have
\begin{align*}
&\xi^{(1)}(s)=\xi^{(1)}(t)-\int_{t}^{s}\nabla_{x}V(\tau,x^{(1)}(\tau))d\tau,
\\&x-y=t\xi^{(1)}(t)-\int_{0}^{t}\int_{t}^{s'}\nabla_{x}V(\tau,x^{(2)}(\tau))d\tau ds'.
\end{align*}
Thus we obtain
\begin{equation} \label{34}
\xi^{(1)}(s)=\frac{x-y}{t}+\frac{1}{t}\int_{0}^{t}\int_{s}^{s'}\nabla_{x}V(\tau,x^{(1)}(\tau))d\tau ds'.
\end{equation}
Similarly, we derive
\begin{align} 
&\xi^{(2)}(s)=\xi-\int_{t}^{s}\nabla_{x}V(\tau,x^{(2)}(\tau))d\tau, \label{35}
\\&
x^{(2)}(0)=x'-t\xi+\int_{0}^{t}\int_{t}^{s'}\nabla_{x}V(\tau,x^{(2)}(\tau))d\tau ds'. \label{136}
\end{align}
It follows from (\ref{136}) that
\begin{align} \label{36}
\frac{x-y}{t}-\xi=\frac{(x-x')-(y-x^{(2)}(0))}{t}-\frac{1}{t}\int_{0}^{t}\int_{t}^{s'}\nabla_{x}V(\tau,x^{(2)}(\tau))d\tau ds'.
\end{align}
Using (\ref{34}), (\ref{35}) and (\ref{36}), we obtain
\begin{align}\label{37}
&\xi^{(1)}(0)-\xi^{(2)}(0)=\frac{x-y}{t}-\xi+\frac{1}{t}\int_{0}^{t}\int_{t}^{s'}\nabla_{x}V(\tau,x^{(1)}(\tau))d\tau ds' \notag
\\
&\hspace{35mm}+\int_{0}^{t}\nabla_{x}V(\tau,x^{(1)}(\tau))d\tau-\int_{0}^{t}\nabla_{x}V(\tau,x^{(2)}(\tau))d\tau \notag
\\
&=\frac{(x-x')-(y-x^{(2)}(0))}{t}+\frac{1}{t}\int_{0}^{t}\int_{t}^{s'}\left\{\nabla_{x}V(\tau,x^{(1)}(\tau))-\nabla_{x}V(\tau,x^{(2)}(\tau))\right\}d\tau ds' \notag
\\
&\hspace{5mm}+\int_{0}^{t}\left\{\nabla_{x}V(\tau,x^{(1)}(\tau))-\nabla_{x}V(\tau,x^{(2)}(\tau))\right\}d\tau \notag
\\
&=\frac{(x-x')-(y-x^{(2)}(0))}{t}+\frac{1}{t}\int_{0}^{t}\int_{0}^{s'}\left\{\nabla_{x}V(\tau,x^{(1)}(\tau))-\nabla_{x}V(\tau,x^{(2)}(\tau))\right\}d\tau ds' \notag
\\
&=: \frac{(x-x')-(y-x^{(2)}(0))}{t}+R_{1}(t,x,y,x',\xi), 
\end{align}
where $R_{1}(t)=R_{1}(t,x,y,x',\xi)$ satisfies, by Taylor's theorem and Assumption \ref{AssV},
\begin{equation*}
|R_1(t)|\lesssim \int_{0}^{t}|x^{(1)}(\tau)-x^{(2)}(\tau)|d\tau.
\end{equation*}
First we prove (\ref{inq1}). We have
\begin{align*}
\frac{d^2}{ds^2}\left(x^{(1)}(s)-x^{(2)}(s)\right)=-\int_{t}^{s}\left(\nabla_{x}V(\tau,x^{(1)}(\tau))-\nabla_{x}V(\tau,x^{(2)}(\tau))\right)d\tau.
\end{align*}
Integrating the both sides of the above equality two times, we obtain
\begin{align*}
x^{(1)}(s)-x^{(2)}(s)=&x^{(1)}(0)-x^{(2)}(0)+s\left(\xi^{(1)}(0)-\xi^{(2)}(0)\right)
\\
&-\int_{0}^{s}\int_{0}^{s'}\int_{t}^{s''}\left\{\nabla_{x}V(\tau,x^{(1)}(\tau))-\nabla_{x}V(\tau,x^{(2)}(\tau))\right\}d\tau ds''ds'.
\end{align*}
This together with (\ref{37}) yields
\begin{align*}
|x^{(1)}(s)-x^{(2)}(s)|&\lesssim |y-x^{(2)}(0)|+t|\xi^{(1)}(0)-\xi^{(2)}(0)|+t^{2}\int_{0}^{t}|x^{(1)}(\tau)-x^{(2)}(\tau)|d\tau
\\&\lesssim |x-x'|+2|y-x^{(2)}(0)|+t|R_1(t)|+t^{2}\int_{0}^{t}|x^{(1)}(\tau)-x^{(2)}(\tau)|d\tau
\\&\lesssim \left(|x-x'|+|y-x^{(2)}(0)|\right)+\int_{0}^{t}|x^{(1)}(\tau)-x^{(2)}(\tau)|d\tau.
\end{align*}
By Gronwall's inequality, we have
\begin{align*}
|x^{(1)}(s)-x^{(2)}(s)|\lesssim |x-x'|+|y-x^{(2)}(0)|,
\end{align*}
which implies
\begin{equation} \label{R1}
|R_1(t)|\lesssim t \left(|x-x'|+|y-x^{(2)}(0)|\right).
\end{equation}
Next we show (\ref{inq2}).
We write $X=(x-x')-(y-x^{(2)}(0))$ for short. 
By the fundamental theorem of calculus and (\ref{37}), we obtain
\begin{align*}
\xi^{(1)}(s)-\xi^{(2)}(s)&=\xi^{(1)}(0)-\xi^{(2)}(0)-\int_{0}^{s}\left\{\nabla_{x}V(\tau,x^{(1)}(\tau))-\nabla_{x}V(\tau,x^{(2)}(\tau))\right\}d\tau
\\&=\frac{X}{t}+R_{1}(t)-\int_{0}^{s}\left\{\nabla_{x}V(\tau,x^{(1)}(\tau))-\nabla_{x}V(\tau,x^{(2)}(\tau))\right\}d\tau.
\end{align*}
It follows from (\ref{R1}) that
\begin{equation*}
\left|\xi^{(1)}(s)-\xi^{(2)}(s)-\frac{X}{t}\right| \lesssim t \left(|x-x'|+|y-x^{(2)}(0)|\right).
\end{equation*}
\end{proof}
%
\subsection{Key estimate}
In this subsection, we show the following lemma which plays an important role in the proof of Theorem \ref{main thm}.
\begin{Lem}\label{Lem 2.5}
For $f \in \CAL{S}(\mathbb{R}^n)$ and multi-indices $\alpha$ with $|\alpha|\leq 2$, there exist positive constants $C_1,C_2$ independent of $t, \e>0$ such that 
\begin{align}
&\left\|x^{\alpha} e^{\frac{1}{2}it\Delta}D_{\varepsilon}f\right\|_{L^1}\leq C_1 \left(\varepsilon+\frac{t}{\varepsilon}\right)^{|\alpha|}\left(\frac{t}{\varepsilon}\right)^{\frac{n}{2}}\hspace{2mm}\sum_{k=0}^{n+1}\left(\frac{\varepsilon}{\sqrt{t}}\right)^{k}, \label{inq10}
\\&
\left\|x^{\alpha} e^{\frac{1}{2}it\Delta}D_{\varepsilon}f\right\|_{L^2}\leq C_2 \left(\e+\frac{t}{\e}\right)^{|\alpha|}. \label{inq11}
\end{align}
\end{Lem}
\begin{proof}
First we prove (\ref{inq10}) in the case $\alpha=0$. We have by (\ref{9}),
\begin{align*} \label{22}
e^{\frac{1}{2}it\Delta}D_{\varepsilon}f(x)&=\int_{\R^n}e^{ix\xi}e^{-\frac{1}{2}it|\xi|^2}\varepsilon^{\frac{n}{2}}\widehat{f}(\varepsilon \xi)d\xi  \notag
\\&=\varepsilon^{\frac{n}{2}}e^{\frac{ix^2}{2t}}\int e^{-\frac{1}{2}it(\xi-\frac{x}{t})^2}\widehat{f}(\varepsilon \xi)d\xi  \notag
\\&=\left(\frac{\varepsilon}{t}\right)^{\frac{n}{2}}e^{\frac{ix^2}{2t}}\int e^{-\frac{1}{2}i|\xi|^2}\widehat{f}\left(\frac{\varepsilon}{\sqrt{t}}\xi+\frac{\varepsilon}{t}x\right)d\xi. 
\end{align*}
Using the equality $e^{-\frac{1}{2}i|\xi|^2}=\displaystyle\frac{1+i\xi \cdot \nabla_{\xi}}{1+|\xi|^2}e^{-\frac{1}{2}i|\xi|^2}$ and integration by parts for $n+1$ times, we have
\begin{align*}
\left\|e^{\frac{1}{2}it\Delta}D_{\varepsilon}f(x)\right\|_{L_{x}^{1}}&\leq \left(\frac{t}{\varepsilon}\right)^{\frac{n}{2}}\sum_{|\beta|\leq n+1}\left(\frac{\varepsilon}{\sqrt{t}}\right)^{|\beta|} \left\|\partial^{\beta}\widehat{f}\right\|_{L^1}\int \langle \xi \rangle^{-n-1}d\xi
\\&\lesssim \left(\frac{t}{\varepsilon}\right)^{\frac{n}{2}}\sum_{|\beta|\leq n+1}\left(\frac{\varepsilon}{\sqrt{t}}\right)^{|\beta|} \left\|\partial^{\beta}\widehat{f}\right\|_{L^1}.
\end{align*}
Next we show (\ref{inq10}) in the case $|\alpha|=1, 2$. Thanks to the facts that
\begin{equation} \label{eq26}
x_{j}^{k}e^{\frac{1}{2}it\Delta}=e^{\frac{1}{2}it\Delta}(x_{j}-it\partial_{x_j})^k,  \hspace{2mm}k=1,2
\end{equation}
and
\begin{align}\label{eq16}
\begin{cases}
&\left(x_{j}-it\partial_{x_j}\right)D_{\varepsilon}=\varepsilon D_{\varepsilon}x_{j}-\displaystyle\frac{it}{\varepsilon}D_{\varepsilon}\partial_{x_j}
\\&\left(x_j-it\partial_{x_j}\right)^{2}D_{\varepsilon}=\varepsilon^2D_{\varepsilon}x_{j}^{2}-itD_{\varepsilon}x_{j}\partial_{x_j}-itD_{\varepsilon}-itD_{\varepsilon}x_j{}-\left(\displaystyle\frac{t}{\varepsilon}\right)^2D_{\varepsilon}\partial_{x_j}^{2},
\end{cases}
\end{align}
we can compute similarly to the case $\alpha=0$.  
Thus we have the desired result with respect to (\ref{inq10}).
Finally (\ref{eq26}), (\ref{eq16}) and the equality $\|e^{\frac{1}{2}it\Delta}D_{\e}f\|_{L^2}=\|f\|_{L^2}$
imply (\ref{inq11}).
\end{proof}
%
%
%
%
\section{Proof of Theorem \ref{main thm}}
In this section, we prove Theorem \ref{main thm} by the following three steps.
%
%
%
%
\subsection{Transformation of equation via wave packet transform}
In this subsection, we transform the Schr\"{o}dinger equation (\ref{CP1}) into a integral equation by using a method similar to that applied to \cite{KKI}.
Let the initial value $u_0\in {\CAL S}({\R}^n)$ and $u(t,x)\in C(\R; L^2(\R^n))$ be the solution to (\ref{CP1}). 
In addition, let $\varphi\in{\CAL S}({\R}^n)\setminus \{0\}$ and $\varphi_{\e}^{(t)}=e^{\frac{1}{2}it\Delta} D_{\e}\varphi$, $t, \e>0$. Taking the wave packet transform $W_{\varphi_{\e}^{(t)}}$ of both sides of (\ref{CP1}), we have
\begin{align}\label{CP1'}
W_{\varphi_{\e}^{(t)}}u(t,x,\xi)&=e^{-i\int_{0}^{t}h(s;t,x,\xi)ds}W_{\varphi_{\e}^{(0)}}u_0(x^{(2)}(0;t,x,\xi),\xi^{(2)}(0;t,x,\xi)) \notag
\\&\hspace{4mm}-i\int_{0}^{t}e^{-i\int_{\tau}^{t}h(s;t,x,\xi)ds}R[u]\left(\tau,x^{(2)}(\tau;t,x,\xi),\xi^{(2)}(\tau;t,x,\xi);\e\right)d\tau,
\end{align}
where $h(s;t,x,\xi)=\frac{1}{2}|\xi^{(2)}(s;t,x,\xi)|^2+V(s,x^{(2)}(s;t,x,\xi))-\nabla_{x}V(s,x^{(2)}(s;t,x,\xi))\cdot x^{(2)}(s;t,x,\xi)$ 
and 
\begin{align}
R[u](t,x,\xi;\e)&=\frac{1}{2}\displaystyle\sum_{|\alpha|=2}\int_{\R^n}\int_{0}^{1}\left(\partial_{x}^{\alpha}V\right)(s,x+\theta(y-x) (1-\theta)d\theta \notag
\\&\hspace{30mm}\times (y-x)^{\alpha}\overline{\varphi_{\e}^{(t)}(y-x)}u(t,y)e^{-iy\xi}dy
\end{align}
(see \cite{KKI}, p. 740-741).
 Taking $W_{\varphi_{\e}^{(t)}}^{*}$ to the both sides of (\ref{CP1'}), we have by Proposition \ref{Prop} and $\|\varphi_{\e}^{(t)}\|_{L^2}=\|\varphi\|_{L^2}$,
\begin{align*}
&u(t,x)=\|\varphi\|_{L^2}^{-2}\iiint_{\R^{3n}} \varphi_{\e}^{(t)}(x-x') \overline{\varphi_{\e}^{(0)}(y-x^{(2)}(0;t,x',\xi))} \notag
\\&
\hspace{40mm}\times e^{-i\int_{0}^{t}h(s;t,x',\xi)ds} e^{-iy\xi^{(2)}(0;t,x',\xi)} e^{ix\xi}u_0(y) dydx'\bar{d}\xi \notag 
\\&
-i\|\varphi\|_{L^2}^{-2}W_{\varphi_{\e}^{(t)}}^{*}\left[\int_{0}^{t}e^{-i\int_{\tau}^{t}h(s;t,x',\xi)ds}R[u]\left(\tau,x^{(2)}(\tau;t,x',\xi),\xi^{(2)}(\tau;t,x',\xi);\e \right)d\tau \right](t,x).
\end{align*}
Concerning first term on the right hand of the above equality, we have by Fubini's theorem,
\begin{align} 
&u(t,x)=\int_{\R^n}E_0(t,x,y;\varepsilon)u_0(y)dy \notag
\\
&-i\|\varphi\|_{L^2}^{-2}W_{\varphi_{\e}^{(t)}}^{*}\left[\int_{0}^{t}e^{-i\int_{\tau}^{t}h(s;t,x',\xi)ds}R[u]\left(\tau,x^{(2)}(\tau;t,x',\xi),\xi^{(2)}(\tau;t,x',\xi);\e \right)d\tau \right](t,x), \label{eq1}
\end{align}
where 
\begin{align}\label{E_0}
&E_0(t,x,y;\varepsilon)=\|\varphi\|_{L^2}^{-2}\iint_{\R^{2n}} \varphi_{\e}^{(t)}(x-x') \overline{\varphi_{\e}^{(0)}(y-x^{(2)}(0;t,x',\xi))} \notag
\\&
\hspace{50mm}\times e^{-i\int_{0}^{t}h(s;t,x',\xi)ds} e^{-iy\xi^{(2)}(0;t,x',\xi)} e^{ix\xi}dx'\bar{d}\xi,
\end{align}
for any $\e>0$.
%
%
\subsection{Estimate for principal part}
In this subsection, our goal is to prove Lemma \ref{Lem 3.1} below.
We choose small $T>0$ satisfing Lemma \ref{Lem1} and \ref{Lem2}.
\begin{Lem} \label{Lem 3.1}
Under Assumption \ref{AssV}, $E_0(t,x,y;\varepsilon)$ is expressed as follows:
\begin{equation*}
E_0(t,x,y;\varepsilon)=(2\pi it)^{-\frac{n}{2}} e^{iS(t,x,y)} \widetilde{a_0}(t,x,y;\varepsilon),\hspace{2mm} (t,x,y)\in(0,T)\times \R^{n+n}, \hspace{2mm}\e>0,
\end{equation*}
where the amplitude $\widetilde{a_0}(t,x,y;\varepsilon)$ satisfies the following estimate:
there exists $C>0$ independent of $t$ and $\e$ such that 
\begin{equation*}
\left\|\widetilde{a_0}(t,x,y;\varepsilon)-1\right\|_{L_{x,y}^{\infty}}\leq C \left(\e+\frac{t}{\e}\right)^2 \hspace{2mm} \sum_{k=0}^{n+1}\left(\frac{\varepsilon}{\sqrt{t}}\right)^{k},\hspace{3mm}t\in(0,T), \e>0.
\end{equation*}
\end{Lem}
we prepare the following lemma in order to prove Lemma \ref{Lem 3.1}.
Let $\phi$ be the phase part of the integrand function for (\ref{E_0}). Namely we put
\begin{equation*}
\phi(t,x,y,x',\xi)=-\int_{0}^{t}h(s;t,x',\xi)ds-y\xi^{(2)}(0;t,x',\xi)+x\xi
\end{equation*}
\begin{Lem}\label{Lem 3.2}
The phase $\phi$ is represented as follows:
\begin{equation*}
\phi(t,x,y,x',\xi)=S(t,x,y)-\frac{1}{2t}\{(x-x')-(y-x^{(2)}(0;t,x',\xi))\}^2+t r_1(t,x,y,x',\xi),
\end{equation*}
where $S(t,x,y)$ is the action integral and $r_0$ satisfyies for some $C>0$,
\begin{equation*} 
|r_1(t,x,y,x',\xi)|\leq C\{(x-x')^2+(y-x^{(2)}(0;t,x',\xi))^2\}.
\end{equation*}
\end{Lem}
\begin{proof}
From the Hamiltonian system (\ref{H}) and integration by parts, we have
\begin{align*}
\int_{0}^{t}x^{(2)}(s)\cdot\nabla_{x}V(s,x^{(2)}(s))ds&=-\int_{0}^{t}x^{(2)}(s)\cdot \frac{d}{ds}\xi^{(2)}(s)ds
\\&=-x'\xi+x^{(2)}(0)\xi^{(2)}(0)+\int_{0}^{t}|\xi^{(2)}(s)|^{2}ds.
\end{align*}
It follows that
\begin{equation*}
\phi(t,x,y,x',\xi)=\int_{0}^{t}\left(\frac{1}{2}|\xi^{(2)}(s)|^{2}-V(s,x^{(2)}(s))\right)ds+(x-x')\xi-(y-x^{(2)}(0))\xi^{(2)}(0).
\end{equation*}
From using equalities
\begin{align*}
&|\xi^{(2)}(s)|^2=|\xi^{(1)}(s)|^2+2\xi^{(1)}(s)\cdot(\xi^{(2)}(s)-\xi^{(1)}(s))+|\xi^{(2)}(s)-\xi^{(1)}(s)|^2,
\\&
V(s,x^{(2)}(s))=V(s,x^{(1)}(s))+\nabla_{x}V(s,x^{(1)}(s))\cdot(x^{(2)}(s)-x^{(1)}(s))+R_{2}(s)
\end{align*}
and employing the equality by integration by parts:
\begin{align*}
&\hspace{5mm}\int_{0}^{t}\nabla_{x}V(s,x^{(1)}(s)(x^{(2)}(s)-x^{(1)}(s))ds
\\&=-\xi^{(1)}(t)(x'-x)+\xi^{(1)}(0)(x^{(2)}(0)-y)+\int_{0}^{t}\xi^{(1)}(s)\cdot(\xi^{(2)}(s)-\xi^{(1)}(s))ds,
\end{align*}
we obtain
\begin{align*}
\phi(t,x,y,x',\xi)&=S(t,x,y)-(x-x')(\xi^{(1)}(t)-\xi)+(y-x^{(2)}(0))(\xi^{(1)}(0)-\xi^{(2)}(0)) \notag
\\&\hspace{5mm}+\frac{1}{2}\int_{0}^{t}|\xi^{(1)}(s)-\xi^{(2)}(s)|^{2}ds-\int_{0}^{t}R_{2}(s)ds,
\end{align*}
where 
\begin{align*}
R_{2}(s)=\frac{1}{2}\displaystyle\sum_{|\alpha|=2}\int_{0}^{1} \partial_{x}^{\alpha}V \left(s,x^{(1)}(s)+\theta(x^{(2)}(s)-x^{(1)}(s))\right) (1-\theta)d\theta \left(x^{(2)}(s)-x^{(1)}(s)\right)^{\alpha}.
\end{align*}
We write $X=(x-x')-(y-x^{(2)}(0))$ for simplicity. From Lemma \ref{Lem2}, we have 
\begin{equation} \label{R2}
|R_2(s)|\lesssim |x-x'|^2+|y-x^{(2)}(0)|^2
\end{equation}
and we can write
\begin{equation} \label{R3}
\xi^{(1)}(s)-\xi^{(2)}(s)=\frac{X}{t}+R_{3}(s)
\end{equation}
with $|R_3(s)|\lesssim t(|x-x'|+|y-x^{(2)}(0)|)$ for all $0\leq s\leq t$.
By (\ref{R3}), we have
\begin{align*}
\phi&=S(t,x,y)-\frac{X^2}{t}-(x-x')R_3(t)+(y-x^{(2)}(0))R_3(0)
\\&
\hspace{5mm}+\frac{1}{2}\int_{0}^{t}\left|\frac{X}{t}+R_3(s)\right|^2ds-\int_{0}^{t}R_2(s)ds
\\&
=S(t,x,y)-\frac{X^2}{2t}
\\&\hspace{5mm}+\left(-(x-x')R_3(t)+(y-x^{(2)}(0))R_3(0)+\int_{0}^{t}\frac{X}{t}R_3(s)+\frac{1}{2}|R_3(s)|^2-R_2(s)ds\right)
\\&=:S(t,x,y)-\frac{X^2}{2t}+t r_1(t,x,y,x',\xi).
\end{align*}
It follows from (\ref{R2}) and (\ref{R3}) that
\begin{equation} \label{r1}
|r_1|\lesssim \left(|x-x'|^2+|y-x^{(2)}(0)|^2\right).
\end{equation}
\end{proof}
\begin{proof}[Proof of Lemma \ref{Lem 3.1}]
From Lemma \ref{Lem 3.2}, we can write 
\begin{equation*}
E_0(t,x,y;\varepsilon)=(2\pi it)^{-\frac{n}{2}}e^{iS(t,x,y)} \widetilde{a_0}(t,x,y;\varepsilon),
\end{equation*}
where 
\begin{align*}
\widetilde{a_0}(t,x,y;\e)=&\|\varphi\|_{L^2}^{-2}\left(\displaystyle\frac{it}{2\pi}\right)^{n/2}\iint_{\R^{2n}} \varphi_{\e}^{(t)}(x-x') \overline{\varphi_{\e}^{(0)}(y-x^{(2)}(0))} \notag
\\&\hspace{20mm}\times e^{-i\frac{1}{2t}\{(x-x')-(y-x^{(2)}(0))\}^2}e^{it r_1(t,x,y,x',\xi)}dx'd\xi.
\end{align*}
By Lemma \ref{Lem1} and the equality $e^{itr_1}=1+(e^{itr_1}-1)$,
we divide the amplitude $\widetilde{a_0}$ into three parts:
\begin{align*}
&\widetilde{a_0}(t,x,y;\varepsilon)=T_1+T_2+T_3,
\\&
T_1=\|\varphi\|_{L^2}^{-2}\left(\frac{it}{2\pi}\right)^{n/2}t^{-n}\iint \varphi_{\e}^{(t)}(x-x')\overline{\varphi_{\e}^{(0)}(y-X)}e^{-i\frac{1}{2t}\{(x-x')-(y-X)\}^2}dXdx',
\\&
T_2=\|\varphi\|_{L^2}^{-2}\left(\frac{it}{2\pi}\right)^{n/2}t^{-n}\iint \varphi_{\e}^{(t)}(x-x')\overline{\varphi_{\e}^{(0)}(y-X)}e^{-i\frac{1}{2t}\{(x-x')-(y-X)\}^2}\left(e^{itr_1}-1\right)dXdx',
\\&
T_3=\|\varphi\|_{L^2}^{-2}\left(\frac{it}{2\pi}\right)^{n/2}t^{-n+1}\iint r_0(t,x',X)\varphi_{\e}^{(t)}(x-x')\overline{\varphi_{\e}^{(0)}(y-X)}e^{-i\frac{1}{2t}\{(x-x')-(y-X)\}^2}dXdx'.
\end{align*}
We obtain by (\ref{9}),
\begin{equation*}
T_1=\|\varphi\|_{L^2}^{-2}\int \varphi_{\e}^{(t)}(x') \left( e^{-\frac{1}{2}it\Delta}\hspace{2mm} \overline{\varphi_{\e}^{(0)}}\right)(x')dx'=\|\varphi\|_{L^2}^{-2}\|\varphi_{\e}^{(t)}\|_{L^2}^{2}=1.
\end{equation*}
It turns out from the inequality $|e^{itr_1}-1|\leq t|r_1|$, (\ref{r1}) and Lemma \ref{Lem 2.5} that
\begin{align*}
|T_2|&\lesssim t^{-n/2+1}\iint \left(|x-x'|^2+|y-X|^2\right)\left|\varphi_{\e}^{(t)}(x-x')\overline{\varphi_{\e}^{(0)}(y-X)}\right|dXdx'
\\&
= t^{-n/2+1}\left(\left\|(x')^{2}\varphi_{\e}^{(t)}(x')\right\|_{L^1}\left\|\varphi_{\e}^{(0)}\right\|_{L^1}+\left\|\varphi_{\e}^{(t)}\right\|_{L^1} \left\|X^{2}\varphi_{\e}^{(0)}(X)\right\|_{L^1}\right)
\\&
\lesssim t\left(\e+\frac{t}{\e}\right)^2 \hspace{2mm}\sum_{k=0}^{n+1}\left(\frac{\varepsilon}{\sqrt{t}}\right)^{k}.
\end{align*}
By Lemma \ref{Lem 2.5}, we have
\begin{align*}
|T_3|\lesssim t^{-n/2+1}\|\varphi_{\e}^{(t)}\|_{L^1} \|\varphi_{\e}^{(0)}\|_{L^1}\lesssim t \sum_{k=0}^{n+1}\left(\frac{\varepsilon}{\sqrt{t}}\right)^{k}.
\end{align*}
Since $\widetilde{a_0}-1=T_2+T_3$, we have desired result.
\end{proof}
%
%
%
%
\subsection{Estimate for remainder term} 
We choose small $T>0$ satisfing Lemma \ref{Lem1} and \ref{Lem2}.
Let $U(t)$ be the propagator for (\ref{CP1}) and  $E(t,x,y)$ be the fundamental solution.
We write $\widetilde{r}(t,x,y;\e)=E(t,x,y)-E_{0}(t,x,y;\e)$.
Then $\widetilde{r}(t,x,y;\e)$ satisfies estimate below.
\begin{Lem} \label{Lem r}
Under Assumption \ref{AssV},
for all $f\in L^2(\R^n), t\in(0,T)$, there exists $C>0$ independent of $f, t$ and $\e$ such that 
\begin{equation*}
\left\| \int_{\R^n} \widetilde{r}(t,\cdot,y;\e)f(y)dy\right\|_{L^{2}(\R^n)} \leq C t\left(\e+\frac{t}{\e}\right)^2  \|f\|_{L^{2}(\R^n)}.
\end{equation*}
\end{Lem}
\begin{proof}
First of all, by density, it is sufficient to consider $f\in \CAL{S}(\R^n)$.
Thanks to (\ref{eq1}), we have for all $f\in \CAL{S}(\R^n)$, 
\begin{align*}
&\left\|\int \widetilde{r}(t,x,y;\e)f(y) dy\right\|_{L^2}
=\left\|U(t)f-\int E_{0}(t,x,y;\e)f(y)dy\right\|_{L^2}
\\&\lesssim \int_{0}^{t}\left\|W_{\varphi_{\e}^{(t)}}^{*}\left[e^{-i\int_{\tau}^{t}h(s;t,x',\xi)ds}R[U(\tau)f](\tau,x^{(2)}(\tau;t,x',\xi),\xi^{(2)}(\tau;t,x',\xi);\e)\right](\tau,t,x)\right\|_{L^2}d\tau.
\end{align*}
From Proposition \ref{Prop}, we continue the above estimate 
\begin{align}\label{inq49}
&\left\|\int \widetilde{r}(t,x,y;\e)f(y) dy\right\|_{L^2} \notag
\\&\lesssim \int_{0}^{t}\|\varphi\|_{L^2}\|R[U(\tau)f](\tau,x^{(2)}(\tau;t,x',\xi),\xi^{(2)}(\tau;t,x',\xi))\|_{L_{x',\xi}^{2}}d\tau  \notag
\\&\lesssim \sum_{|\alpha|=2}\int_{0}^{t}\left\|\int_{{\mathbb R}^n}\int_{0}^{1} \left(\partial_{x}^{\alpha}V\right)\left(\tau,x^{(2)}(\tau;t,x',\xi)+\theta(y-x^{(2)}(\tau;t,x',\xi))\right)(1-\theta)d\theta  \right.\notag
\\& \left. \hspace{5mm}\times (y-x^{(2)}(\tau;t,x',\xi))^{\alpha}\overline{\varphi_{\e}^{(\tau)}(y-x^{(2)}(\tau;t,x',\xi))}U(\tau)f(y)e^{-iy \xi^{(2)}(\tau;t,x',\xi)}dy\right\|_{L_{x',\xi}^{2}}d\tau. 
\end{align}
We now consider the change of variables $X'=x^{(2)}(\tau;t,x',\xi)$ and $\varXi=\xi^{(2)}(\tau;t,x',\xi)$. 
Since the Jacobian for the change of variables $(x',\xi)\to(X',\Xi)$ equals $1$, it follows from the Plancherel's theorem and Assumption \ref{AssV} that
\begin{align}\label{inq50}
&\left\|\int_{{\mathbb R}^n}\int_{0}^{1} \left(\partial_{x}^{\alpha}V\right)\left(\tau,x^{(2)}(\tau;t,x',\xi)+\theta(y-x^{(2)}(\tau;t,x',\xi))\right)(1-\theta)d\theta  \right. \notag
\\& \left. \hspace{5mm}\times (y-x^{(2)}(\tau;t,x',\xi))^{\alpha}\overline{\varphi_{\e}^{(\tau)}(y-x^{(2)}(\tau;t,x',\xi))}(U(\tau)f)(y)e^{-iy \xi^{(2)}(\tau;t,x',\xi)}dy\right\|_{L_{x',\xi}^{2}}d\tau. \notag
\\&=\left\|\int_{{\mathbb R}^n} \int_{0}^{1}\left(\partial_{x}^{\alpha} V\right)(\tau,X'+\theta(y-X')) (1-\theta)d\theta (y-X')^{\alpha}\overline{\varphi_{\e}^{(\tau)}(y-X')}(U(\tau)f)(y)e^{-iy \varXi}dy\right\|_{L_{X',\Xi}^{2}} \notag
\\&=\left\|\int_{0}^{1} \left(\partial_{x}^{\alpha} V\right)(\tau,X'+\theta(\varXi-X')) (1-\theta)d\theta (\varXi-X')^{\alpha}\overline{\varphi_{\e}^{(\tau)}(\varXi-X')}(U(\tau)f)(\varXi)\right\|_{L_{X',\varXi}^{2}}  \notag
\\&\lesssim \left\|(\varXi-X')^{\alpha}\overline{\varphi_{\e}^{(\tau)}(\varXi-X')}(U(\tau)f)(\varXi)\right\|_{L_{X',\varXi}^{2}}  \notag
\\&=\left\|x^{\alpha} \varphi_{\e}^{(\tau)}(x)\right\|_{L_{x}^{2}} \left\|U(\tau)f\right\|_{L^2}.  
\end{align}
Plugging (\ref{inq50}) into (\ref{inq49}), we estimate, by conservation law of $U(\tau)$ in $L^2$ and Lemma \ref{Lem 2.5},
\begin{align*}
\left\|\int \widetilde{r}(t,x,y;\e)f(y) dy\right\|_{L^2}
&\lesssim \sum_{|\alpha|=2}\int_{0}^{t} \left(\e+\frac{\tau}{\e}\right)^2 \|U(\tau)f\|_{L^2}d\tau 
\\&\lesssim t\left(\e+\frac{t}{\e}\right)^2 \|f\|_{L^2}.  \notag
\end{align*}
\end{proof}
We define $a_{0}(t,x,y)=\widetilde{a_0}(t,x,y;\sqrt{t})$ and $r(t,x,y)=\widetilde{r}(t,x,y;\sqrt{t})$.
By (\ref{eq1}), we have (\ref{main-form}). Moreover by Lemma \ref{Lem 3.1} and \ref{Lem r}, we obtain the estimates (\ref{inqa_0}) and (\ref{inqr}) of Theorem \ref{main thm}: we complete the proof of Theorem \ref{main thm}.
%

%
%
%
\end{document}